\documentclass[12pt]{amsart}
\renewcommand{\bar}{\overline}
\usepackage{amsmath,amsthm,amscd,euscript,amssymb}
\renewcommand{\bar}{\overline}
\def \r{\mathbb R}

\DeclareMathOperator{\codim}{codim}
\DeclareMathOperator{\sgn}{sgn} 

\newtheorem{theorem}{Theorem}[section]
\newtheorem{lemma}[theorem]{Lemma}

\newtheorem{proposition}[theorem]{Proposition}
\newtheorem{corollary}[theorem]{Corollary}
\theoremstyle{remark}
\newtheorem{remark}[theorem]{Remark}
\theoremstyle{definition}
\newtheorem{definition}[theorem]{Definition}
\newtheorem{example}[theorem]{Example}
\newtheorem{problem}{Problem}
\newtheorem{conjecture}[problem]{Conjecture}

\title{Geometry of configuration spaces of tensegrities}

\author[Franck Doray \and Oleg Karpenkov \and Jan Schepers]{Franck Doray\,$^1$, Oleg Karpenkov\,$^2$ and Jan Schepers\,$^1$}
\thanks{$^1$ Supported by VICI-grant 639.033.402 of NWO}

\thanks{$^2$ Partially supported by RFBR grant SS-709.2008.1 and by NWO-DIAMANT grant 613.009.001.}

\keywords{Tensegrity, self-tensional equilibrium frameworks,
stratification}

\email[Franck Doray]{doray@math.leidenuniv.nl}
\email[Oleg Karpenkov]{karpenk@mccme.ru}
\email[Jan Schepers]{janschepers1@gmail.com}

\address[Franck Doray, Oleg Karpenkov, Jan Schepers]{Mathematisch Instituut, Universiteit Leiden,
P.O. Box 9512, 2300 RA Leiden, The Netherlands}

\begin{document}

\input epsf

\begin{abstract}
Consider a graph $G$ with $n$ vertices. In this paper we study
geometric conditions for an $n$-tuple of points in $\r^d$ to admit
a tensegrity with underlying graph $G$. We introduce and
investigate a natural stratification, depending on $G$, of the
configuration space of all $n$-tuples in $\r^d$. In particular we
find surgeries on graphs that give relations between different
strata. Based on numerous examples we give a description of
geometric conditions defining the strata for plane tensegrities,
we conjecture that the list of such conditions is sufficient to
describe any stratum. We conclude the paper with particular
examples of strata for tensegrities in the plane with a small
number of vertices.
\end{abstract}

\maketitle

\tableofcontents

\section{Introduction}

In his paper~\cite{Max} J.~C.~Maxwell made one of the first
approaches to the study of equilibrium states for frames under
the action of static forces. He noted that the frames together
with the forces give rise to reciprocal figures. In the second
half of the twentieth century the artist K.~Snelson built many
surprising sculptures consisting of cables and bars that are
actually such frames in equilibrium, see~\cite{Sne}.
R.~Buckminster Fuller introduced the name ``tensegrity'' for
these constructions, combining the words ``tension'' and
``integrity''. A nice overview of the history of tensegrity
constructions is made by R.~Motro in his book~\cite{Mot}.

In mathematics, tensegrities were investigated in several papers.
In~\cite{Rot1} B.~Roth and W.~Whiteley and in~\cite{Con}
R.~Connelly and W.~Whiteley studied rigidity and flexibility of
tensegrities, see also the survey about rigidity in~\cite{Whi}.

N.~L.~White and W.~Whiteley started in~\cite{WW} the investigation
of geometric realizability conditions for a tensegrity with
prescribed bars and cables. In the preprint~\cite{Guz2}
M.~de~Guzm\'an describes other examples of geometric conditions
for tensegrities.

Tensegrities have a wide range of applications in different
branches of science and in architecture. For instance they are
used in the study of viruses~\cite{Cas}, cells~\cite{Ing}, for
construction of deployable mechanisms~\cite{Ske,Tib}, etc.

\vspace{2mm}

We focus on the following important question. Suppose a graph $G$
is given. {\it Is the graph $G$ realizable as a tensegrity for a
general configuration of its vertices?} We develop a new
technique to study this question. We introduce special operations
({\it surgeries}) that change the graph in a certain way but
preserve the property to be (not to be) realizable as a
tensegrity.

Let $n$ be the number of vertices of $G$. Consider the
configuration space of all $n$-tuples of points in $\r^d$. In
this paper we define a stratification of the configuration space
such, that each stratum corresponds to a certain set of admissible
tensegrities associated to $G$. Suppose that one wants to obtain
a construction with some edges of $G$ replaced by struts and the
others by cables, then he/she should take a configuration in a
specific stratum of the stratification.

In this paper we prove that all the strata are semialgebraic
sets, and therefore a notion of dimension is well-defined for
them. This allows to generalize the previous question: {\it what
is the minimal codimension of the strata in the configuration
space that contains $n$-tuples of points admitting a tensegrity
with underlying graph $G$?} Our technique of surgeries on graphs
also gives the first answers in this case. In particular we
obtain the list of all 6, 7, and 8 vertex tensegrities in the
plane that are realizable for codimension~1 strata. We note that
the complete answers to the above questions are not known to the
authors.

N.~L.~White and W.~Whiteley~\cite{WW} and M.~de~Guzm\'an and
D.~Orden~\cite{Guz3,Guz1} have found the geometric conditions of
realizability of plane tensegrities with $6$ vertices and of some
other particular cases. We continue the investigation for other
graphs (see Subsection~6.2). In all the observed examples the
strata are defined by certain systems of geometric conditions. It
turns out that all these geometric conditions are obtained from
elementary ones:

--- {\it two points coincide};

--- {\it three points are on a line};

--- five points $a$, $b$, $c$, $d$, $e$ satisfy: {\it $e$ is the
intersection point of the lines passing through points $a$ and
$b$ and
points $c$ and $d$ respectively}.\\
We conjecture that for plane tensegrities any stratum can be
defined by certain geometric conditions (see Section~5).

\vspace{2mm}

This paper is organized as follows. We start in Section~2 with
general definitions. In Subsection~2.1 we describe the
configuration space of tensegrities associated to a given graph as
a fibration over the affine space of all frameworks. We introduce
a natural stratification on the space of all frameworks in
Subsection~2.2. We prove that all strata are semialgebraic sets
and therefore the strata have well-defined dimensions. In
Section~3 we study the dimension of solutions for graphs on
general configurations of points in $\r^d$. Later in this section
we calculate the dimensions in the simplest cases, and formulate
general open questions. In Section~4 we study surgeries on graphs
and frameworks that induce isomorphisms of the spaces of
self-stresses for the frameworks. We give general definitions
related to systems of geometric conditions for plane tensegrities
in Section~5. We conjecture that any stratum is a dense subset of
the solution of one of such systems. Finally in Section~6 we give
particular examples of graphs and their strata for tensegrities
in the plane. We study the dimension of the space of
self-stresses in Subsection~6.1 and give tables of geometric
conditions for codimension 1 strata for graphs with 8 vertices
and less in Subsection~6.2.

\vspace{2mm}

{\bf Acknowledgements.} The authors are grateful to B.~Edixhoven
for rousing our interest to the subject, to A.~Sossinski,
V.~Goryunov, and A.~Perucca for helpful remarks and discussions,
to S.~Speed for useful information on graph classification, and
Mathematisch Instituut of Universiteit Leiden for hospitality and
excellent working conditions. The second author is grateful to
Liverpool University for the organization of a fruitful visit.

\section{General definitions}

\subsection{Configuration spaces of tensegrities}

Recall a slightly modified definition of a framework
from~\cite{Guz1}.

\begin{definition}
Fix a positive integer $d$. Let $G=(V,E)$ be an arbitrary graph
without loops and multiple edges. Let it have $n$ vertices.
\itemize
\item{A {\it framework} $G(P)$ in $\r ^d$ is a map of the graph $G$ with vertices $v_1, \ldots, v_n$ on
a finite point configuration $P=(p_1,\ldots,p_n)$ in $\r^d$ with
straight edges, such that $G(P)(v_i)=p_i$ for $i=1,\ldots,n$.}
\item{A {\it stress} $w$ on a framework is an assignment of real scalars
$w_{i,j}$ (called {\it tensions}) to its edges $p_ip_j$. We also
put $w_{i,j}=0$ if there is no edge between the corresponding
vertices. Observe that $w_{i,j}=w_{j,i}$, since they refer to the
same edge.}

\item{A stress $w$ is called a {\it self-stress} if, in addition, the
following equilibrium condition is fulfilled at every vertex
$p_i$:
$$
\sum\limits_{\{j|j\ne i\}} w_{i,j}\bar{p_ip_j}=0.
$$
By $\bar{p_ip_j}$ we denote the vector from the point $p_i$ to the
point $p_j$.}
\item{A couple $(G(P),w)$ is called a {\it tensegrity} if
$w$ is a self-stress  for the framework $G(P)$.}
\end{definition}

\begin{remark} Tensegrities are self-tensional equilibrium frameworks. For
instance, any framework for the two vertex graph without edges is
always a tensegrity, although it is not rigid. For more
information about rigidity of tensegrities we refer to a paper of
B.~Roth and W.~Whiteley~\cite{Rot1}.
\end{remark}

Denote by $W(n)$ the linear space of dimension $n^2$ of all
tensions $w_{i,j}$. Consider a framework $G(P)$ and denote by
$W(G,P)$ the subset of $W(n)$ of all possible self-stressed
tensions for $G(P)$. By definition of self-stressed tensions, the
set $W(G,P)$ is a linear subspace of $W(n)$.

The {\it configuration space of tensegrities} corresponding to
the graph $G$ is the set
$$
\big\{(G(P), w)\,|\,P\in (\r^d)^n, w \in W(G,P)\big\},
$$
we denote it by $\Omega(G)$. The set $\{G(P)\,|\,P\in (\r^d)^n\}$
is said to be the {\it base of the configuration space}, we
denote it by $B_d(G)$. If we forget about the edges between the
points in all the frameworks, then we get natural bijections
between $\Omega(G)$ and a subset of $(\r^d)^n\times W(n)$ and
between $B_d(G)$ and $(\r^d)^n$. Later on we actually identify
the last two pairs of sets. The bijections induce natural
topologies on $\Omega(G)$ and $B_d(G)$.

Let $\pi$ be the natural projection of $\Omega(G)$ to the base
$B_d(G)$. This defines the structure of a fibration. For a given
framework $G(P)$ of the base we call the set $W(G,P)$ {\it the
linear fiber} at the point $P$ (or at the framework $G(P)$) of the
configuration space.

Consider a self-stress $w$ for the framework $G(P)$. We say that
the edge $p_ip_j$ is a {\it cable} if $w_{i,j}<0$ and a {\it
strut} if $w_{i,j}>0$.

\begin{remark}\label{light}
The definitions of struts and cables come from the following
physical interpretation. Suppose we would like to construct the
lightest possible tensegrity structure on a given framework and
with a given self-stress using heavy struts and relatively light
cables. Then we should replace the edges with positive $w_{i,j}$
with struts, and the edges with negative $w_{i,j}$ with cables.
Such constructions would be the lightest possible.
\end{remark}

Denote by ``$\sgn$'' the sign function over $\r$.

\begin{definition}
Consider a framework $G(P)$ and one of its self-stresses $w$. The
$n\times n$ matrix $(\sgn(w_{i,j}))$ is called the {\it
strut-cable matrix} of the stress $w$ and denoted by $\sgn(w)$.
\end{definition}

Let us give one example of a strut-cable matrix.

\begin{example}
Consider a configuration of four points in the plane: $p_1(0,0)$,
$p_2(1,0)$, $p_3(2,2)$, $p_4(0,1)$ and a self-stress $w$ as on the
picture: $w_{1,2}=6$, $w_{1,3}=-3$, $w_{1,4}=6$, $w_{2,3}=2$,
$w_{2,4}=-4$, $w_{3,4}=2$. Then we have:
$$
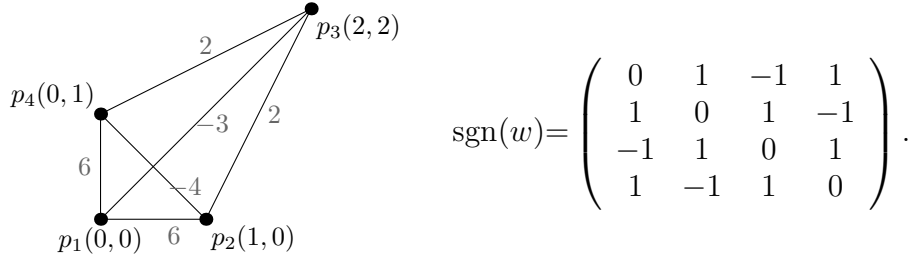
\begin{array}{c}\epsfbox{graphs.1}\end{array}
\quad
\sgn(w){=}\left(
\begin{array}{cccc}
 0& 1&-1& 1\\
 1& 0& 1&-1\\
-1& 1& 0& 1\\
 1&-1& 1& 0\\
\end{array}
\right) .
$$
\end{example}

\subsection{Stratification of the base of a configuration space of tensegrities}

Suppose we have some framework $G(P)$ and we want to find the
lightest cable-strut construction on it, as explained in
Remark~\ref{light}. Then the following questions arise. {\it
Which edges can be replaced by cables, and which by struts?} {\it
What is the geometric position of points in the configurations
for which given edges may be replaced by cables and the others by
struts?} The questions lead to the following definition.

\begin{definition}
A linear fiber $W(G,P_1)$ is said to be {\it equivalent} to a
linear fiber $W(G,P_2)$ if there exists a homeomorphism $\xi$
between $W(G,P_1)$ and $W(G,P_2)$, such that for any self-stress
$w$ in $W(G,P_1)$ the self-stress $\xi(w)$ satisfies
$$
\sgn\big(\xi(w)\big)=\sgn\big(w\big).
$$
\end{definition}

The described equivalence relation gives us a stratification of
the base $B_d(G)=(\r^d)^n$. A {\it stratum} is by definition a
maximal connected set of points with equivalent linear fibers.
Once we have proven Theorem~\ref{semialg}, by general theory of
semialgebraic sets (see for instance~\cite{Boc}) it follows that
all strata are path-connected.

\begin{example}
We describe the stratification of $B_1(K_3)=\r^3$ for the
complete graph $K_3$ on three vertices. The point $(x_1,x_2,x_3)$
in $\r^3$ corresponds to the framework with vertices $p_1=(x_1)$,
$p_2=(x_2)$, and $p_3=(x_3)$. The stratification consists of 13
strata. There is 1 one-dimensional stratum, and there are 6
two-dimensional and 6 three-dimensional strata.

The one-dimensional stratum consists of frameworks with all
vertices coinciding. It is defined by the equations
$x_1=x_2=x_3$. The dimension of the fiber at a point of this
stratum is 3.

Any of the two-dimensional strata consists of frameworks with
exactly two vertices coinciding. The strata are the connected
components of the complement to the line $x_1=x_2=x_3$ in the
union of the three planes $x_1=x_2$, $x_1=x_3$, and $x_2=x_3$.
The dimension of the fiber at a point of any of these strata is 2.

Any of the three-dimensional strata consists of frameworks with
distinct vertices. The strata are the connected components of the
complement in $\r^3$ to the union of the three planes $x_1=x_2$,
$x_1=x_3$, and $x_2=x_3$. The dimension of the fiber at a point
of any of these strata is 1.
\end{example}

In general we have the following theorem.

\begin{theorem}\label{semialg}
Any stratum is a semialgebraic set.
\end{theorem}

For the definition and basic properties of semialgebraic sets we
refer the reader to~\cite{Boc}.

\vspace{1mm}

We need two preliminary lemmas for the proof of the theorem, but
first we introduce the following notation.

Let $M$ be an arbitrary symmetric $n{\times}n$-matrix with zeroes
on the diagonal and all the other entries belonging to
$\{-1,0,1\}$. Let $i$ be an integer with $0\le i\le n^2$. We say
that a couple $(M,i)$ is a {\it stratum symbol}.

For an arbitrary framework $G(P)$ we denote by $W_M(G,P)$ the set
of all self-stresses with strut-cable matrix $\sgn(w)$ equal to
$M$. The closure of $W_M(G,P)$ is a pointed polyhedral cone with
vertex at the origin. The set $W_M(G,P)$ is homeomorphic to an
open $k$-dimensional disc, we call $k$ the dimension of $W_M(G,P)$
and denote it by $\dim(W_M(G,P))$.

For any stratum symbol $(M,i)$ we denote by $\Xi(M,i)$ the set
$$
\{(G(P),w)\,|\,w\in W(G,P), \sgn(w)=M,\dim(W_M(G,P))=i\} \subset
\Omega(G).
$$

\begin{lemma}\label{ll1}
For any stratum symbol $(M,i)$, the subset $\pi(\Xi(M,i))$ of the
base $B_d(G)$ is either empty or it is a semialgebraic set.
\end{lemma}

\begin{proof}
The set $\Xi(M,i)$ is a semialgebraic set since it is defined by
a system of equations and inequalities in the coordinates of
the vertices and the tensions of the following three types:\\
{\it a$)$} quadratic equilibrium condition equations;\\
{\it b$)$} linear equations or inequalities specifying if
the coordinate values $w_{i,j}$ are zeroes, positive, or negative reals;\\
{\it c$)$} algebraic equations and inequalities defining
respectively $\dim(W_M(G,P))\le i$ and $\dim(W_M(G,P))\ge i$.
Note that $\dim(W_M(G,P))$ is equal to the dimension of the
linear space spanned by $W_M(G,P)$.

Let us make a small remark about item {\it $($c$)$}. At each
framework we take the system of equilibrium conditions and
equations of type $w_{i,j}=0$ in the variables $w_{i,j}$. This
system consists of the {\it equalities} of items {\it $($a$)$}
and {\it $($b$)$}. It is linear in the variables $w_{i,j}$. The
coefficients of the corresponding matrix depend linearly on the
coordinates of the framework vertices. The equations and
inequalities of item {\it $($c$)$} are defined by some
determinants of submatrices of this matrix being equal or not
equal to zero. Therefore, they are algebraic.

Since by the Tarski-Seidenberg theorem any projection of a
semialgebraic set is semialgebraic, the set $\pi(\Xi(M,i))$ is
semialgebraic.
\end{proof}

Denote by $S(G,P)$ the set of all stratum symbols $(M,i)$ that are
realized by the point $G(P)$, in other words
$$
S(G,P)=\{(M,i)\,|\,G(P)\in \pi(\Xi(M,i))\}.
$$

\begin{lemma}\label{ll2}
Let $G(P_1)$ and $G(P_2)$ be two frameworks. Then
$S(G,P_1)=S(G,P_2)$ if and only if the linear fiber $W(G,P_1)$ is
equivalent to the linear fiber $W(G,P_2)$.
\end{lemma}

\begin{proof}
Let the linear fiber at the point $G(P_1)$ be equivalent to the
linear fiber at the point $G(P_2)$ then by definition we have
$$
S(G,P_1)=S(G,P_2).
$$

\vspace{1mm}

Suppose now that $S(G,P_1)=S(G,P_2)$. Let us denote by
$\bar{W(G,P_i)}$ the one point compactification of the fiber
$W(G,P_i)$ for $i=1,2$. So $\bar{W(G,P_i)}$ is homeomorphic to a
sphere of dimension $\dim W(G,P_i)$.

For any point $P$ and any $M$ the set $W_M(G,P)$ is a convex cone
homeomorphic to an open disc of dimension $\dim(W_M(G,P))$. So,
for any point $P$ we have a natural CW-decomposition of
$\bar{W(G,P)}$ with cells $W_M(G,P)$ and the new one point cell.

A cell $W_{M'}(G,P_1)$ is adjacent to a cell $W_{M''}(G,P_1)$ iff
the cell $W_{M'}(G,P_2)$ is adjacent to the cell
$W_{M''}(G,P_2)$. This is true, since the couples of cells
corresponding to $M'$ and to $M''$ are defined by the same sets
of equations and inequalities of type ``$>$'', and the closures
of $W_{M'}(G,P_i)$  for $i=1,2$ are defined by the system
defining $W_{M'}(G,P_i)$ with all ``$>$'' in the inequalities
replaced by ``$\ge$''.

Therefore, there exists a homeomorphism of $\bar{W(G,P_1)}$ and
$\bar{W(G,P_2)}$, sending all the cells $W_{M}(G,P_1)$ to the
corresponding cells $W_{M}(G,P_2)$. We leave the proof of this
statement as an exercise for the reader, this can be done by
inductively constructing the homeomorphism on the $k$-skeletons
of the CW-complexes.

Hence, the linear fibers $W(G,P_1)$ and $W(G,P_2)$ are equivalent.
\end{proof}

\vspace{1mm}

{\it Proof of Theorem~\ref{semialg}.} Let us prove the theorem for
a stratum containing some point~$P$. Consider any point $P'$ in
the stratum. By definition, $W(G,P)$ is equivalent to the space
$W(G,P')$, and hence by Lemma~\ref{ll2}, we have $S(G,P)=S(G,P')$.

Consider the following set
$$
\bigcap\limits_{(M,i)\in S(G,P)}\pi(\Xi(M,i))\setminus
\Bigl(\bigcup\limits_{(M,i)\notin S(G,P)}\pi(\Xi(M,i))\Bigr),
$$
we denote it $\Sigma(P)$. So $\Sigma(P)$ is the set of frameworks
$G(P')$ for which $S(G,P')=S(G,P)$. By Lemma~\ref{ll1} all the
sets $\pi(\Xi(M,i))$ are semialgebraic. Therefore, the set
$\Sigma(P)$ is semialgebraic. Denote by $\Sigma'(P)$ the connected
component of $\Sigma(P)$ that contains the point $P$. Since the
set $\Sigma(P)$ is semialgebraic, the set $\Sigma'(P)$ is also
semialgebraic, see~\cite{Boc}.

Let us show that $\Sigma'(P)$ is the stratum of $B_d(G)$
containing the point $P$. First, the set $\Sigma'(P)$ is
contained in the stratum. This holds since $\Sigma'(P)$ is
connected and consists of points with equivalent sets $S(G,P)$.
And hence by Lemma~\ref{ll2} all the points of $\Sigma'(P)$ have
equivalent linear fibers $W(G,P)$. Secondly, the stratum is
contained in the space $\Sigma'(P)$. This holds since the stratum
is connected and consists of points with equivalent linear fibers
$W(G,P)$. Thus by Lemma~\ref{ll2} all the points of the stratum
have equivalent sets $S(G,P)$.

As we have shown, the stratum containing $P$ coincides with
$\Sigma'(P)$ and hence it is semialgebraic. \qed

\vspace{2mm}

From the above proof it follows that the total number of strata
is finite.

\section{On the tensegrity $d$-characteristic of graphs}

In this section we study the dimension of the linear fiber for
graphs on a general point configuration in $\r^d$. We give a
natural definition of the tensegrity $d$-characteristic of a
graph and calculate it for the simplest graphs. In addition we
formulate general open questions for further investigation.

\subsection{Definition and basic properties of the tensegrity $d$-characteristic}

Note that for any two points $P_1$ and $P_2$ of the same stratum
$S$ of the space $B_d(G)$ for a graph $G$ we have
$$\dim(W(G,P_1))=\dim(W(G,P_2)).$$
Denote this number by $\dim(G,S)$. Denote also by $\codim(S)$ the
integer
$$\dim(B_d(G)){-}\dim(S).$$
Consider a graph $G$ with at least one edge. We call the integer
$$
\min\{\codim S\,|\,\hbox{$S$ is a stratum of $B_d(G)$, }
\dim(G,S)>0\}
$$
the {\it codimension} of $G$ and denote it by $\codim_d(G)$.

\begin{definition}
We call the integer
$$
\left\{
\begin{array}{l}
1-\codim_d(G)  \hbox{, if $\codim_d(G)>0$}
\\
\max\big\{\dim(W(G,P))\,\big|\, \hbox{$G(P)$ contained in a
codimension zero stratum}\big\}, \hbox{ otherwise}
\end{array}
\right.
$$
the {\it tensegrity $d$-characteristic} of the graph $G$ (or the
{\it $d$-TC} of $G$ for short), and denote it by $\tau_d(G)$.
\end{definition}

%
%

\begin{example}
Consider the two graphs shown on Figure~\ref{graphs.2}. The left
one is a graph of codimension~1 in the plane, it can be realized
as a tensegrity iff either the two triangles are in perspective
position or the points of one of the two triples $(v_1,v_4,v_5)$
or $(v_2,v_3,v_6)$ lie on a line (for more details
see~\cite{Guz1}), so its 2-TC is zero. The graph on the right has
a twodimensional space of self-stresses for a general position
plane framework, and hence its 2-TC equals two (we show this
later in Proposition~\ref{Proposition 2.2.}).

\begin{figure}
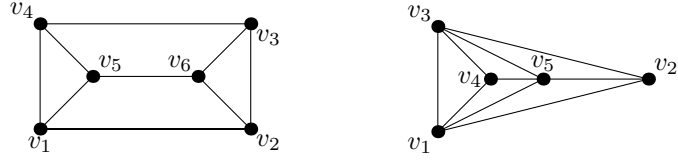

$$\epsfbox{graphs.2} \qquad \qquad \epsfbox{graphs.3}$$
\caption{A graph with zero tensegrity 2-characteristic (the left
one) and a graph whose 2-characteristic equals 2 (the right
one).}\label{graphs.2}
\end{figure}

\end{example}

\begin{proposition}
Let $S_1$ and $S_2$ be two strata of codimension 0. Let $G(P_1)$
and $G(P_2)$ be two points of the strata $S_1$ and $S_2$
respectively. Then the following holds:
$$
\dim(W(G,P_1))=\dim(W(G,P_2)).
$$
\end{proposition}

\begin{proof}
The equilibrium conditions give a linear system of equations in
the variables $w_{i,j}$, at each framework linearly depending on
the coordinates of the vertices. The dimension of the solution
space is determined by the rank of the matrix of this system. The
subset of $B_d(G)$ where the rank is not maximal is an algebraic
subset of positive codimension. By definition, this set does not
have elements in the strata $S_1$ and $S_2$. This yields the
statement of the proposition.
\end{proof}

\begin{corollary}\label{codimzerostrata}
Let $G$ be a graph.
If $\tau_d(G)\geq 0$ then for every framework $G(P)$ in a
codimension 0 stratum we have $\dim W(G,P)=\tau_d(G)$. \qed
\end{corollary}


\subsection{Atoms and atom decomposition}

In this subsection we recall a definition and some results of
M.~de~Guzm\'an and D.~Orden~\cite{Guz1} that we use later.

Consider a point configuration $P$ of $d{+}2$ points in general
position in $\mathbb{R}^{d}$. Throughout this subsection `general
position' means that no $d{+}1$ of them are contained in a
hyperplane. An {\it atom} in $\mathbb{R}^d$ is a tensegrity
$(K_{d+2}(P),w)$, where $K_{d+2}$ is the complete graph on
$d{+}2$ vertices and where $w$ is a nonzero self-stress.

According to \cite[Section 2]{Guz1} the linear fiber
$W(K_{d+2},P)$ is one-dimensional for $P$ in general position, in
particular this implies $\tau_d(K_{d+2})=1$. In addition the
tension on every edge in the atom is nonzero. A more general
statement holds.

\begin{lemma}\label{Lemma1.2.}
\cite[Lemma 2.2]{Guz1} Let $G(P)$ be a framework on a point
configuration $P$ in $\mathbb{R}^d$ in general position. Let
$p\in P$. Given a nonzero self-stress on $G(P)$, then either at
least $d{+}1$ of the edges incident to $p$ receive nonzero
tension, or all of them have zero tension.
\end{lemma}

M.~de~Guzm\'an and D.~Orden showed that one can consider atoms as
the building blocks of tensegrity structures. First, we explain
how to add tensegrities. Let $T=(G(P),w)$ and $T'=(G'(P'),w')$ be
two tensegrities. We define $T+T'$ as follows. The framework of
$T+T'$ is $G(P)\cup G'(P')$, we take the union of vertices and
edges. The tension on a common edge $p_ip_j=p'_kp'_l$ is defined
as $w_{i,j}+w'_{k,l}$ and on an edge appearing exactly in one of
the original frameworks we put the original tension. It is easy to
see that the defined stress is a self-stress, so $T+T'$ is a
tensegrity.

\begin{theorem}
\cite[Theorem 3.2]{Guz1} Every tensegrity $(G(P),w)$ with  a
general position point configuration $P$ and $w_{i,j}\neq 0$ on
all edges of $G$ is a finite sum of atoms. This decomposition is
not unique in general.
\end{theorem}

\subsection{Calculation of tensegrity $d$-characteristic in the simplest cases}

We start this subsection with the formulation of a problem, we do
not know the complete solution of it.

\begin{problem}
Give a general formula for $\tau_d(G)$ in terms of the
combinatorics of the graph.
\end{problem}

Let us calculate the $d$-TC for a complete graph, this will give
us the maximal value of the $d$-TC for fixed number of vertices
$n$ and dimension $d$.

\begin{proposition}\label{Proposition 1.5.}
For any positive integers $n$ and $d$ satisfying $n\geq d{+}2$,
we have
$$
\tau_d(K_n) = \frac{(n-d-1)(n-d)}{2}.
$$
\end{proposition}

\begin{proof} We work by induction on $n$. For $n=d+2$ the $d$-TC equals
1, as mentioned above. For $n>d{+}2$ we choose any point
configuration $P$ on $n$ points such that no $d{+}1$ of them lie
in a hyperplane. Take $p\in P$. Any tensegrity $(K_n(P),w)$ can be
decomposed as a sum of $n{-}d{-}1$ atoms with $p$ as vertex and a
tensegrity on $P\setminus\{p\}$ with underlying graph $K_{n-1}$.
Indeed, we can use such atoms to cancel the given tensions on
$n{-}d{-}1$ edges at $p$. Then there are only $d$ edges left, so
by Lemma~\ref{Lemma1.2.} the tensions on these edges equal zero.
We conclude by induction that

$$
\tau_d(K_n)  = \tau_d(K_{n-1}) + n-d-1 = \frac{(n-d-1)(n-d)}{2}.
$$

\end{proof}

Now we show how the $d$-TC behaves when we remove an edge of the
graph.

\begin{proposition}\label{0_-1}
Let $G$ be some graph satisfying $\tau_d(G)>1$. Let a graph $G'$
be obtained from the graph $G$ by erasing one edge. Then
$$
\tau_d(G)-\tau_d(G') \in \{0, 1\}.
$$
\end{proposition}

\begin{proof}
Erasing one edge is equivalent to adding a new linear equation
$w_{i,j}=0$ to the linear system defining the space $W(G,P)$ for
the graph $G$ (for any point $P$). This implies that the space of
solutions coincides with $W(G,P)$ or it is a hyperplane in
$W(G,P)$.

So, first, $\tau_d(G')\le\tau_d(G)$.

Secondly, since $\tau_d(G)=\dim(W(G,P_0))$ for some framework
$P_0$ of a codimension 0 stratum for $G$ (and therefore it
belongs to a codimension 0 stratum for $G'$), then
$$
\tau_d(G')=\dim(W(G',P_0))\ge\dim(W(G,P_0))-1=\tau_d(G)-1.
$$

This completes the proof.
\end{proof}

As we show in the example below, erasing an edge does not always
reduce the tensegrity characteristic.

\begin{example}
Consider the graph shown in Figure~\ref{1}. Assume that this graph
underlies a tensegrity. Then we can add an atom on the four
leftmost vertices, to cancel the tension on edge $f$ for
instance. This automatically cancels the tensions on the edges
connecting the four leftmost vertices by Lemma~\ref{Lemma1.2.}.
We can do the same on the right. So the tension on $e$ is zero as
well. Therefore the tension on $e$ was zero from the beginning and
hence deleting $e$ does not change the $2$-TC. In
Example~\ref{Example 2.3.} we give a less trivial example of this
phenomenon.
\begin{figure}
$$
\epsfbox{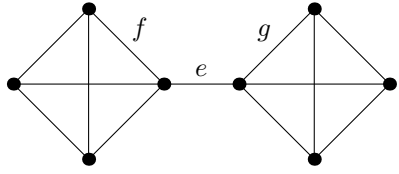}
$$
\caption{The tension on the edge $e$ is always zero.}\label{1}
\end{figure}
\end{example}

Let us formulate two general corollaries of
Proposition~\ref{0_-1}.

\begin{corollary}\label{Corollary 1.7.}
Let $G$ be a graph on $n$ vertices and $m\in \mathbb{Z}_{>0}$. If
$G$ has
$$
m + \frac{n(n-1)}{2} - \tau_d(K_n) = m + dn -\frac{d^2+d}{2}
$$
edges, then $\tau_d(G) \geq m$.
\end{corollary}

\begin{proof}
Combine Proposition~\ref{Proposition 1.5.} and
Proposition~\ref{0_-1}.
\end{proof}

The following corollary is useful for the calculation of the
tensegrity $d$-characteristic. In Subsection~\ref{smallvert1} we
use it to calculate all the tensegrity $2$-characteristics for
sufficiently connected graphs with less than 8 vertices.

\begin{corollary}\label{Corollary 1.8.}
Let $G$ be a graph on $n$ vertices with $\tau_d(G)\geq 0$. Assume
that $G$ has
$$
dn -\frac{d^2+d}{2}+\tau_d(G)
$$
edges. Then for any
graph $H$ that can be obtained from $G$ by adding $N$ edges we
have
$$
\tau_d(H)= \tau_d(G) + N.
$$
\end{corollary}

\begin{proof}
We delete  $\tau_d(K_n)-\tau_d(G)-N$ edges from $K_n$ to reach
$H$. If the $d$-TC does not drop by 1 at one of these steps, then
we apply Proposition~\ref{0_-1} an additional $N$ times to $H$ to
reach $G$. This leads to a wrong value of $\tau_d(G)$. So the
$d$-TC drops by one in each of the first
$\tau_d(K_n)-\tau_d(G)-N$ steps and the formula for $\tau_d(H)$
follows.
\end{proof}

\begin{example}\label{pseudo} A {\it pseudo-triangle} is a
planar polygon with exactly three vertices at which the angles
are less than $\pi$. Let $G$ be a planar graph with $n$ vertices
and $k$ edges that admits a pseudo-triangular embedding $G(P)$ in
the plane, i.e.\ a non-crossing embedding such that the outer
face is convex and all interior faces are pseudo-triangles.  It is
obvious that a pseudo-triangular embedding $G(P)$ belongs to a
codimension 0 stratum of $B_2(G)$. By Lemma~2 of~\cite{Ord1} we
find that

--- $\tau_2(G) = k-(2n-3)$ if $k-(2n-3)\geq 1$,

--- $\tau_2(G) \leq 0$ if $k-(2n-3)=0$.
\\
(Note that for pseudo-triangular embeddings we always have $k\geq
2n{-}3$.)
\end{example}

\section{Surgeries on graphs that preserve the dimension of the fibers}

In this section we describe operations that one can perform on a
graph without changing the dimensions of the corresponding fibers
for the frameworks. We refer to such operations as {\it
surgeries}. The first type of surgeries is for general dimension,
while the other two are restricted to dimension $d=2$. We do not
know other similar operations that are not compositions of the
surgeries described below.

The idea of surgeries is analogous to the idea of Reidemeister
moves in knot theory. If two graphs are connected by a sequence of
surgeries, then one obtains tensegrities for the first graph from
tensegrities for the second graph and vice versa.

We essentially use surgeries to calculate the list of geometric
conditions for the strata for (sufficiently connected) graphs
with less than~9 vertices and with zero $2$-TC in
Subsection~\ref{smallvert2}.

\subsection{General surgeries in arbitrary dimension}
For an edge $e$ of a graph $G$ we denote by $G_e$ the graph
obtained from $G$ by removing $e$. Recall that a subgraph $G'$ of
a graph $G$ is said to be {\it induced} if, for any pair of
vertices $v_i$ and $v_j$ of $G'$, $v_iv_j$ is an edge of $G'$ if
and only if $v_iv_j$ is an edge of $G$.

Denote by $\Sigma_d(G)$ the union of codimension zero strata in
$B_d(G)$. Let $G$ be a graph and $H$ a subgraph. Consider the map
that takes a framework for $G$ to the framework for $H$ by
forgetting all the vertices and edges of $G$ that are not in $H$.
Denote by $\Sigma_d(G,H)$ the preimage of $\Sigma_d(H)$ for this
map.

\begin{proposition}\label{Proposition 1.9.}
Let $G$ be a graph and $H$ an induced subgraph with $\tau_d(H)=1$.
Consider a configuration $P_0$ lying in $\Sigma_d(G,H)$. Suppose
that there exists a self-stress on the framework $G(P_0)$ that has
nonzero tensions for all edges of $H$ and zero tensions on the
other edges. Let $e_1,e_2$ be edges of $H$. Then for any $P\in
\Sigma_d(G,H)$ we have
$$
W(G_{e_1},P) \cong W(G_{e_2},P).
$$
\end{proposition}

The corresponding surgery takes the graph $G_{e_1}$ to $G_{e_2}$,
or vice versa.

\begin{remark}
We always have the inclusion $\Sigma_d(G) \subset \Sigma_d(G,H)$,
this follows directly from the definition of the strata.
Nevertheless the set $\Sigma_d(G,H)$ usually contains many strata
of $B_d(G)$ of positive codimension. So
Proposition~\ref{Proposition 1.9.} is applicable to all strata of
codimension zero as well as to some strata of positive
codimension.
\end{remark}

For the proof of Proposition~\ref{Proposition 1.9.} we need the
following lemma.

\begin{lemma}\label{lll1}
Let $G$ be a graph with $\tau_d(G)=1$ and $e$ one of its edges.
Suppose that there exists a configuration $P_0\in \Sigma_d(G)$ and
a nonzero self-stress $w_0$ such that $w_0(e)=0$. Then for any
tensegrity $(G(P),w)$ with $P\in \Sigma_d(G)$ we get $w(e)=0$.
\end{lemma}

\begin{proof}
Since $\tau_d(G)=1$ and $P_0\in \Sigma_d(G)$, any tensegrity
$(G(P_0),w)$ satisfies the condition $w(e)=0$. Therefore, any
tensegrity with $P$ in the same stratum as $P_0$ has zero tension
at $e$. So the condition {\it always to have zero tension at $e$}
defines a somewhere dense subset $S$ in $B_d(G)$. Since the
condition is defined by a solution of a certain linear system,
$S$ is dense in $B_d(G)$. It follows that $\Sigma_d(G)$ is a
subset of $S$.
\end{proof}

\vspace{2mm}

{\it Proof of Proposition~\ref{Proposition 1.9.}.} From
Lemma~\ref{lll1} we have that for any configuration of
$\Sigma_d(H)$ there exist a unique up to a scalar self-stress
that is nonzero at each edge of $H$. The uniqueness follows from
the fact that $\tau_d(H)=1$. Hence for any configuration of
$\Sigma_d(G,H)$ there exists a unique up to a scalar self-stress
that is nonzero at each edge of $H$ and zero at all other edges
of $G$.

For any $P\in \Sigma_d(G,H)$ we obtain an isomorphism between
$W(G_{e_1},P)$ and $W(G_{e_2},P)$ by adding the unique tensegrity
on the underlying subgraph $H$ of $G$ that cancels the tension on
$e_2$, considered as edge of $G_{e_1}$. \qed

\vspace{2mm}

In particular one can use atoms (i.e.~$H=K_{d+2}$) in the above
proposition.

\begin{corollary}\label{Corollarygeneral}
In the notation and with the conditions of
Proposition~\ref{Proposition 1.9.} we have: if either
$\tau_d(G_{e_1})>0$ or $\tau_d(G_{e_2})>0$ then
$$
\tau_d(G_{e_1}) = \tau_d(G_{e_2}).
$$
\end{corollary}

\begin{proof}
The statement follows directly from Proposition~\ref{Proposition
1.9.} and Corollary~\ref{codimzerostrata}.
\end{proof}

Let us show how to use the above corollary to compute the
tensegrity characteristic.

\begin{example}\label{Example 1.10.}
We calculate the 2-TC of the graph $G_0$ shown in Figure~\ref{2}.
Consider the atom $H$ on the vertices $v_1$, $v_2$, $v_3$, and
$v_4$ and let $e_1$, $e_2$ be the edges $v_2v_4$, $v_1v_3$
respectively. Denote by $G$ the graph obtained from $G_0$ by
adding the edge $e_2$. So the graph $G_0$ is actually $G_{e_1}$.
By Corollary~\ref{Corollary 1.7.} we have $\tau_2(G_0)\geq 1$, and
hence it is possible to apply Corollary~\ref{Corollarygeneral}.
Consider the graph $G_{e_2}$, it is shown in Figure~\ref{2} in
the middle. The degree of the vertex $v_1$ in this graph equals 2,
so by Lemma~\ref{Lemma1.2.} the tensions on its incoming edges
equal zero if the points $v_1$, $v_2$, and $v_4$ are not on a
line. After removing these two edges and the vertex $v_1$ we get
the graph of an atom. Therefore,
$$
\tau_2(G_0)=1.
$$
\begin{figure}
$$
\epsfbox{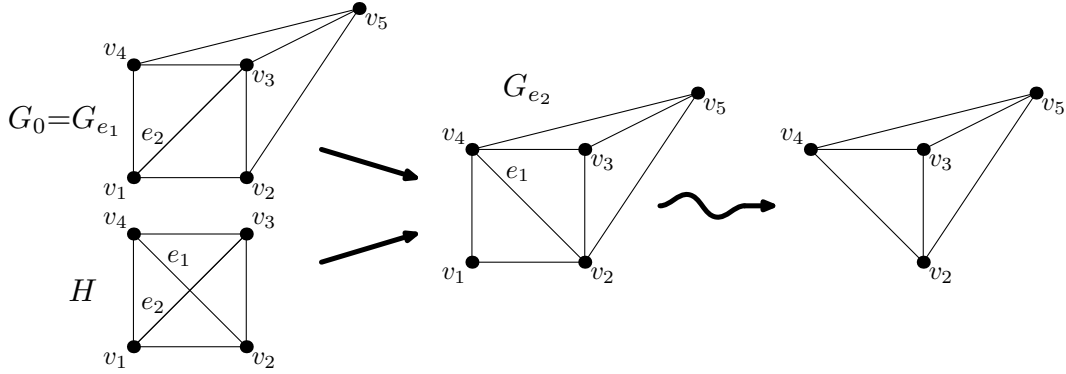}$$ \caption{This shows that
$\tau_2(G_0)=1$.}\label{2}
\end{figure}
\end{example}

\subsection{Additional surgeries in dimension two}\label{oper12}


In this subsection we study two surgeries on edges of plane
frameworks that do not change the dimension of the fibers of the
frameworks.

\vspace{2mm}

{\bf Surgery I.} Consider a graph $G$ and a framework $G(P)$. Let
$G$ contain the complete graph $K_4$ with vertices $v_1$, $v_2$,
$v_3$, and $v_4$ as an induced subgraph. Suppose that the edges
between $v_1,v_2,v_3,v_4$ and other vertices of $G$ are as
follows:

--- $pv_2$ and $qv_3$ for unique vertices $p$ and $q$;

--- the edges $pv_1$ and $qv_1$;

--- any set of edges from $v_4$.

\vspace{2mm}

\begin{figure}
$$\epsfbox{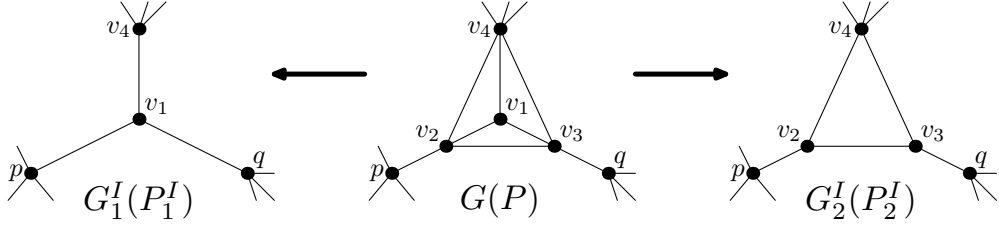}$$
\caption{Surgery I.}\label{operations.1}
\end{figure}

In addition we require that the framework $G(P)$ has the triples
of points $(p,v_1,v_2)$ and $(q,v_1,v_3)$ on one line. See
Figure~\ref{operations.1} in the middle.

Let us delete from the graph $G$ the vertices $v_2$ and $v_3$
(the vertex $v_1$) with all edges adjacent to them. We denote the
resulting graph by $G_1^{I}$ (by $G_2^{I}$ respectively). The
corresponding framework is denoted by $G_1^{I}(P_1^{I})$ (by
$G_2^{I}(P_2^{I})$ respectively). See Figure~\ref{operations.1}
on the left (on the right). Surgery~I takes $G_1^I$ to $G_2^I$ or
vice versa.

\begin{proposition}\label{operation1}
Consider the frameworks $G(P)$, $G_1^{I}(P_1^{I})$, and
$G_2^{I}(P_2^{I})$ as above. If the triples of points
$(p,v_2,v_3)$, $(q,v_2,v_3)$, $(p,v_2,v_4)$, $(q,v_3,v_4)$ and
$(v_2,v_3,v_4)$ are not on a line then we have
$$
W(G_1^{I},P_1^{I}) \cong W(G_2^{I},P_2^{I}).
$$
\end{proposition}

\begin{proof}
We explain how to go from $W(G_2^{I},P_2^{I})$ to
$W(G_1^{I},P_1^{I})$. The inverse map is simply given by the
reverse construction. By the conditions the intersection point
$v_1$ of $pv_2$ and $qv_3$ is uniquely defined and not on the
lines through $v_2$ and $v_4$ or $v_3$ and $v_4$. We add the
uniquely defined atom on $v_1,v_2,v_3,v_4$ to $G_2^{I}(P_2^{I})$
that cancels the tension on $v_2v_3$. Since $p,v_2,v_1$ lie on
one line, this surgery also cancels the tension on $v_2v_4$ and
similarly for $v_3v_4$. Due to the equilibrium condition at
$v_2$, we can replace the edges $pv_2$ and $v_2v_1$ with their
tensions $w_{p,2}$ and $w_{2,1}$ by an edge $pv_{1}$ with tension
$w_{p,1}$ defined by one of the following vector equations:
$$
w_{p,2}\bar{pv_2}=w_{p,1}\bar{pv_1}=w_{2,1}\bar{v_2v_1}.
$$
This uniquely defines a self-stress on $G_1^{I}(P_1^{I})$.
\end{proof}

\begin{corollary}\label{Corollary 1.12.NEW}
Assume that one of the following conditions holds:
\begin{enumerate}
\item $\tau_2(G_1^{I})> 0$ or $\tau_2(G_2^{I})> 0$.
\item $\tau_2(G_1^{I})=0$ and there is a codimension 1 stratum $S$ of $B_2(G_1^I)$ such that

--- $\dim W(G_1^I,P)>0$ for a $G_1^I(P)$ in the stratum $S$,

--- the stratum $S$ is not contained in the subset of $B_2(G_1^I)$ of frameworks
having one of the triples of points $(p,v_1,q)$, $(p,v_1,v_4)$, or
$(q,v_1,v_4)$ on one line.
\item $\tau_2(G_2^{I})=0$ and there is a codimension 1 stratum $S'$ of $B_2(G_2^I)$ such that

--- $\dim W(G_2^I,P')>0$ for a $G_2^I(P')$ in the stratum $S'$,

--- the stratum $S'$ is not contained in the subset of $B_2(G_2^I)$ of
frameworks having $(p,v_2,v_3)$, $(q,v_2,v_3)$, $(p,v_2,v_4)$,
$(q,v_3,v_4)$, or $(v_2,v_3,v_4)$ on one line.

\end{enumerate}
Then
$$
\tau_2(G_1^{I})=\tau_2(G_2^{I}).
$$
\end{corollary}

\begin{proof}
Let $A$ be the subset of $B_2(G_2^I)$ of frameworks having
$(p,v_2,v_3)$, $(q,v_2,v_3)$, $(p,v_2,v_4)$, $(q,v_3,v_4)$ or
$(v_2,v_3,v_4)$ on one line. Let $B$ be the subset of
$B_2(G_1^I)$ of frameworks having $(p,v_1,q)$, $(p,v_1,v_4)$ or
$(q,v_1,v_4)$ on one line. Note that $A$ and $B$ are of
codimension 1. The proof of Proposition~\ref{operation1} gives a
surjective map
$$
\varphi: B_2(G_2^I) \setminus A \to B_2(G_1^I) \setminus B
$$
inducing an isomorphism between the linear fibers above $G(P)\in
B_2(G_2^I) \setminus A$ and $\varphi(G(P))$. Now in all the cases
(1)---(3) the statement of the corollary follows directly from
the definition of the tensegrity characteristic.
\end{proof}


\vspace{2mm}

\begin{figure}
$$\epsfbox{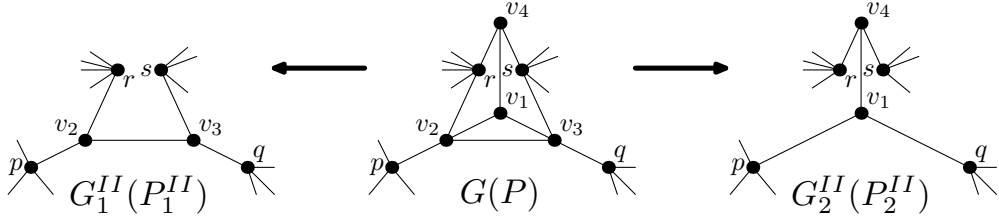}$$
\caption{Surgery II.}\label{operations.2}
\end{figure}

{\bf Surgery II.}  Consider a graph $G$ and a framework $G(P)$.
Let $G$ contain the complete graph $K_4$ with vertices $v_1$,
$v_2$, $v_3$, and $v_4$ as an induced subgraph. Suppose that the
set of edges between $v_1,v_2, v_3,v_4$ and other vertices of $G$
is
$$
\{pv_1, pv_2, qv_1, qv_3, rv_2, rv_4, sv_3, sv_4\},
$$
for unique points $p,q,r,s$. In addition we require that the
framework $G(P)$ has the triples of points
$$
(p,v_1,v_2), \quad (q,v_1,v_3), \quad (r,v_2,v_4), \quad
\hbox{and} \quad (s,v_3,v_4)
$$
on one line. See Figure~\ref{operations.2} in the middle.

Let us delete from the graph $G$ the vertices $v_1$ and $v_4$
($v_2$ and $v_3$) with all edges adjacent to them. We denote the
resulting graph by $G_1^{II}$ (by $G_2^{II}$ respectively). The
corresponding framework is denoted by $G_1^{II}(P_1^{II})$ (by
$G_2^{II}(P_2^{II})$ respectively). See Figure~\ref{operations.2}
on the left (on the right). Surgery~II takes $G_1^{II}$ to
$G_2^{II}$ or vice versa.

The proofs of the proposition and corollary below are similar to
the proofs of Proposition~\ref{operation1} and
Corollary~\ref{Corollary 1.12.NEW}.

\begin{proposition}\label{operation2}
Consider the frameworks $G(P)$, $G_1^{II}(P_1^{II})$, and
$G_2^{II}(P_2^{II})$ as above. If non of the triples of points
$(p,q,v_1)$, $(p,v_1,v_4)$, $(r,v_1,v_4)$, $(q,v_1,v_4)$,
$(s,v_1,v_4)$, or $(r,s,v_4)$ lie on a line then we have
$$
W(G_1^{II},P_1^{II}) \cong W(G_2^{II},P_2^{II}).
$$
\qed
\end{proposition}

\begin{corollary}
Assume that one of the following conditions holds:
\begin{enumerate}
\item $\tau_2(G_1^{II})> 0$ or $\tau_2(G_2^{II})> 0$.
\item $\tau_2(G_1^{II})=0$ and there is a codimension 1 stratum $S$ of $B_2(G_1^{II})$ such that

--- $\dim W(G_1^{II},P)>0$ for a $G_1^{II}(P)$ in the stratum $S$,

--- the stratum $S$ is not contained in the subset of $B_2(G_1^{II})$ of
frameworks having $(p,v_2,v_3)$, $(q,v_2,v_3)$, $(p,v_2,r)$,
$(q,v_3,s)$, $(r,v_2,v_3)$, or $(s,v_2,v_3)$ on one line.

\item $\tau_2(G_2^{II})=0$ and there is a codimension 1 stratum $S'$ of $B_2(G_2^{II})$ such that

--- $\dim W(G_2^{II},P')>0$ for a $G_2^{II}(P')$ in the stratum $S'$,

--- the stratum $S'$ is not contained in the subset of $B_2(G_2^{II})$ of
frameworks having $(p,q,v_1)$, $(p,v_1,v_4)$, $(r,v_1,v_4)$,
$(q,v_1,v_4)$, $(s,v_1,v_4)$, or $(r,s,v_4)$ on one line.
\end{enumerate}
Then
$$
\tau_2(G_1^{II})=\tau_2(G_2^{II}).
$$
\qed
\end{corollary}


\section{Geometric relations for strata and complexity of tensegrities
in two-dimensional case}

In all the observed examples of plane tensegrities with a given
graph the strata for which a tensegrity is realizable are defined
by certain geometric conditions on the points of the
corresponding frameworks. In this section we study such geometric
conditions. In Subsection~\ref{A simple examples} we describe an
example of a geometric condition for a particular graph. Further,
in Subsection~\ref{system} we give general definitions related to
systems of geometric conditions. Finally, in
Subsections~\ref{geometric_structure} and~\ref{complexity} we
formulate two open questions related to the geometric nature of
tensegrity strata.

\vspace{2mm}

To avoid problems with describing annoying cases of
parallel/nonparallel lines we extend the plane $\r^2$ to the
projective space. It is convenient for us to consider the
following model of the projective space: $\r P^2=\r^2\cup
l_\infty$. The set of points $l_\infty$ is the set of all
``directions'' in the plane. The set of lines of $\r P^2$ is the
set of all plane lines (each plane line contains now a new point
of $l_\infty$ that is the direction of $l$) together with the line
$l_\infty$. Now any two lines intersect at exactly one point.

\subsection{A simple example}\label{A simple examples}

First, we study the graph shown in Figure~\ref{graphs.2} on the
left, we denote it by $G_0$. In~\cite{WW} N.~L.~White and
W.~Whiteley proved that the $2$-TC of this graph is zero. They
showed that there exists a nonzero tensegrity with graph $G_0$ and
framework $P$ iff the points of $P$ satisfy one of the following
three
conditions:\\
{\it i$)$} the lines $v_1v_2$, $v_3v_4$, and $v_5v_6$ have a common nonempty intersection (in $\r P^2$);\\
{\it ii$)$} the vertices $v_1$, $v_4$, and $v_5$ are in one line;\\
{\it iii$)$} the vertices $v_2$, $v_3$, and $v_6$ are in one line.

We remind that the base $B(G_0)$ of the configuration space is
$\r^{12}$ with coordinates $(x_1,y_1,\ldots,x_{6},y_6)$, where
$(x_i,y_i)$ are the coordinates of $v_i$. Condition {\it $($i$)$}
defines a degree 4 hypersurface with equation
$$
\det \left(
\begin{array}{ccc}
 y_1-y_2& y_3-y_4& y_5-y_6\\
 x_2-x_1& x_4-x_3& x_6-x_5\\
 x_1y_2-x_2y_1 & x_3y_4-x_4y_3 & x_5y_6-x_6y_5\\
\end{array}
\right)=0.
$$
and Conditions {\it $($ii$)$} and {\it $($iii$)$} define the
conics
$$
\begin{array}{l}
x_1y_4+x_4y_5+x_5y_1-x_1y_5-x_4y_1-x_5y_4=0, \quad \hbox{and}\\
x_2y_3+x_3y_6+x_6y_2-x_2y_6-x_3y_2-x_6y_3=0
\end{array}
$$
respectively.

\subsection{Systems of geometric conditions}\label{system}
%
%

Let us define three elementary geometric conditions. Consider an
ordered subset $P=\{ p_1,\ldots, p_n\}$ of the projective plane.

\vspace{2mm}

\underline{\it 2-point condition}. We say that the subset $P$ {\it
satisfies the condition $p_i=p_j$} if $p_i$ coincides with $p_j$.

\underline{\it 3-point condition}. We say that the subset $P$
{\it satisfies the condition}
$$
p_i\triangledown p_j\triangledown p_k =0
$$
if the points $p_i$, $p_j$, and $p_k$ are on a line.

\underline{\it 5-point condition}. We say that the subset $P$ {\it
satisfies the condition}
$$
p_i=[p_j,p_{j'};p_k,p_{k'}]
$$
if the four points $p_j,p_{j'},p_k$, and $p_{k'}$ are on a line
and $p_i$ also belongs to this line, or if $p_i=p_jp_{j'}\cap
p_kp_{k'}$ otherwise. We say that $[p_j,p_{j'};p_k,p_{k'}]$ is
the {\it intersection symbol} of the lines $p_jp_{j'}$ and
$p_kp_{k'}$.


Note that we define the last condition in terms of closures, since
$[p,q;r,s]$ is not defined for all 4-tuples, but for a dense
subset.

\begin{definition}
Consider a system of elementary geometric conditions for ordered
$n$-point subsets of $\r P^2$, and let $m\le n$.

--- We say that the ordered $n$-point subset $P$ of projective plane
{\it satisfies} the system of elementary geometric conditions if
$P$ satisfies each of these conditions.

--- We say that the ordered subset $\{p_1,\ldots,p_m\}$ {\it
satisfies conditionally} the system of elementary geometric
conditions if there exist points $q_1, \ldots, q_{n-m}$ such that
the ordered set
$$
\{p_1,\ldots,p_m,q_1,\ldots, q_{n-m}\}
$$
satisfies the system. We call the number $n{-}m$ {\it the
conditional number} of the system.
\end{definition}

\begin{example}\label{6pconic}
The condition that {\it six points $p_1, \ldots, p_6$  lie on a
conic} is equivalent to the following geometric conditional
system:
$$
\left\{
\begin{array}{l}
q_1=[p_1,p_{2};p_4,p_5]\\
q_2=[p_2,p_{3};p_5,p_6]\\
q_3=[p_3,p_{4};p_1,p_6]\\
q_1\triangledown q_2\triangledown q_3 =0
\end{array}
\right. .
$$
This is a reformulation of Pascal's theorem. The conditional
number is 3 here.

We can rewrite the system as follows, for short:
$$
[p_1,p_{2};p_4,p_5]\triangledown [p_2,p_{3};p_5,p_6]\triangledown
[p_3,p_{4};p_1,p_6] =0.
$$
\end{example}

\begin{example}\label{2_2_2}
The condition for six points $p_1, \ldots, p_6$ that {\it the
lines $p_1p_2$, $p_3p_4$, and $p_5p_6$ have a common point} is
equivalent to the following geometric conditional system:
$$
\left\{
\begin{array}{l}
q_1=[p_1,p_{2};p_3,p_4]\\
q_1\triangledown p_5\triangledown p_6 =0
\end{array}
\right. ,
$$
or in a shorter form:
$$
[p_1,p_{2};p_3,p_4]\triangledown p_5\triangledown p_6 =0.
$$
The conditional number of the system is 1.
\end{example}

\subsection{Conjecture on geometric structure of the
strata}\label{geometric_structure}

For a given positive integer $k$ and a graph $G$ consider the set
of all frameworks $G(P)$ at which the dimension of the fiber
$W(G,P)$ is greater than or equal to $k$. We call this set the
$(G,k)$-{\it stratum}. Since any $(G,k)$-stratum is a finite union
of strata of the base $B_2(G)$, it is semialgebraic.

\begin{definition}
Let $G$ be a graph and $k$ be a positive integer. The
$(G,k)$-stratum is said to be {\it geometric} if it is a finite
union of the sets of conditional solutions of systems of
geometric conditions (in these systems $p_1,\ldots,p_m$
correspond to the vertices of the graph).
\end{definition}

\begin{conjecture}
For any graph $G$ and integer $k$ the $(G,k)$-stratum is
geometric.
\end{conjecture}

The conjecture is checked for all the graphs with seven and fewer
vertices, see Section~6 for the techniques.

\begin{problem}
Find analogous elementary geometric conditions in the three-
(higher-) dimensional case.
\end{problem}

We refer to~\cite{WW} for examples of geometric conditions in
dimension 3.

\subsection{Complexity of the strata}\label{complexity}

We end this section with a discussion of the complexity of
geometric $(G,k)$-strata.

A geometric $(G,k)$-stratum is defined by some union of the
conditional solutions of systems of geometric conditions. Each
system in this union has its own conditional number. Take the
maximal among all the conditional numbers in the union. We call
the minimal number among such maximal numbers for all the unions
of systems defining the same $(G,k)$-stratum the {\it geometric
complexity} of the $(G,k)$-stratum.

\begin{example}
The geometric complexity of $(G_0,1)$ stratum for the graph $G_0$
described in Subsection~\ref{A simple examples} and shown in
Figure~\ref{graphs.2} on the left equals $3$.
\end{example}

\begin{problem}
Find the asymptotics of the maximal complexity of geometric
$(G,k)$-strata with bounded number of vertices $k$ while $k$
tends to infinity.
\end{problem}


\section{Plane tensegrities with a small number of vertices}

In this section we work in the two-dimensional case (unless
otherwise stated). In Subsection~\ref{smallvert1} we study the
$2$-TC of graphs. In particular, we calculate the $2$-TC for
sufficiently connected graphs with seven or less vertices. In
Subsection~\ref{smallvert2} we give a list of geometric
conditions for realizability of tensegrities in the plane for
graphs with zero $2$-TC.

\subsection{On the tensegrity $2$-characteristic of graphs}\label{smallvert1}

Recall the following definitions from graph theory. Let $G$ be a
graph. The {\it vertex connectivity} $\kappa(G)$ is the minimal
number of vertices whose deletion disconnects $G$. The {\it edge
connectivity} $\lambda(G)$ is the minimal number of edges whose
deletion disconnects $G$.
It is well known that $\kappa(G)\leq \lambda(G)$.

\vspace{2mm}

For general dimension $d$, let $G(P)$ be a framework in
$\mathbb{R}^d$ with underlying graph $G$. If $\kappa(G) <d$ or
$\lambda(G)<d+1$ then $G(P)$ consists of two or more pieces that
can rotate with respect to each other. So for us the most
interesting graphs are those with $\kappa(G) \geq d$ and
$\lambda(G)\geq d{+}1$.

\begin{proposition}\label{Proposition 2.2.}
Let $G$ be a 2-vertex and 3-edge connected graph with $k$ edges
and $n$ vertices. If $n\leq 7$, then
$$
\tau_2(G) = k-2n+3.
$$
\end{proposition}

\begin{remark}
In particular we have equality in Corollary~\ref{Corollary 1.7.}
under the conditions of Proposition~\ref{Proposition 2.2.}. The
formula of Proposition~\ref{Proposition 2.2.} holds for many
graphs in general, see for instance Example~\ref{pseudo}. It does
not always hold for graphs with 9 vertices as the example below
shows.
\end{remark}

\begin{example}\label{Example 2.3.}
Let $G$ be the graph with 9 vertices and 15 edges as in
Figure~\ref{9}. If we use the formula of
Proposition~\ref{Proposition 2.2.}, then we have $\tau_2(G)=0$.
Nevertheless, $G$ contains $K_4$ as an induced subgraph. Hence for
any framework $G(P)$ the dimension of $W(G,P)$ is at least 1 (we
put zero tensions on all edges not belonging to $K_4$ and choose
a nonzero self-stress on $K_4$). So $\tau_2(G)\geq 1$.  In fact
it is not hard to prove that $\tau_2(G)=1$. This in particular
implies that the tensions on all edges not belonging to $K_4$ are
zero for a framework in a codimension zero stratum.

\begin{figure}
$$\epsfbox{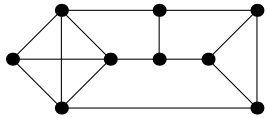}$$
\caption{A graph $G$ with 9 vertices, 15 edges and $\tau_2(G)=
1$.}\label{9}
\end{figure}
\end{example}

Notice that the graph $G$ of Example~\ref{Example 2.3.} is not a
{\it Laman graph}, i.e.\ a graph with $2n{-}3$ edges, where $n$
is the number of vertices, for which each subset of $m\geq 2$
vertices spans at most $2m{-}3$ edges. Theorem 1.1 of \cite{Rot2}
shows that every planar Laman graph $H$ can be embedded as a
pseudo-triangulation and hence $\tau_2(H)\leq 0$ by
Example~\ref{pseudo}. We suspect that equality holds here, and
more generally for all Laman graphs.

\vspace{2mm}

%
%
%
%
%

{\it Proof of Proposition~\ref{Proposition 2.2.}.} We use a
classification argument.

\vspace{1mm}

{\bf Four vertices.} For the complete graph $K_4$ we have
$\tau_2(K_4)=1=6{-}8{+}3$. There are no other graphs satisfying
the conditions of the proposition.

\vspace{1mm}

\begin{figure}
$$\epsfbox{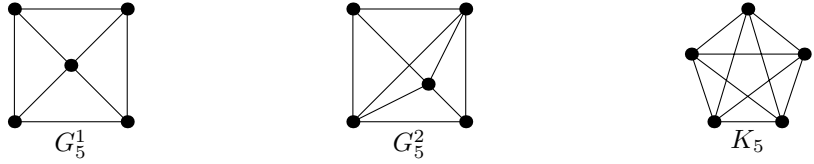}$$
\caption{The three possible graphs with five vertices, $\kappa\geq
2$ and $\lambda\geq 3$.}\label{4}
\end{figure}

{\bf Five vertices.} There are three possibilities, we show them
in Figure~\ref{4}. From Proposition~\ref{Proposition 1.5.}  we
know that $\tau_2(K_5)=3=10-10+3$ and in Example~\ref{Example
1.10.} we have seen that $\tau_2(G_5^1)=1=8-10+3$. To see that
$\tau_2(G_5^2)=2$ we apply Corollary~\ref{Corollary 1.8.}.

\vspace{1mm}

\begin{figure}
$$\epsfbox{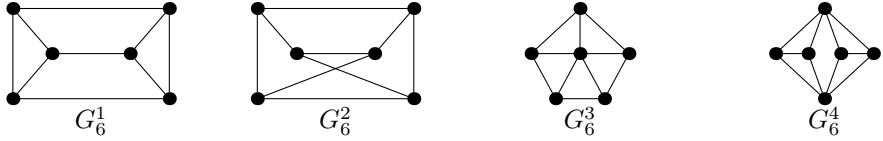}$$
\caption{The four graphs with six vertices, $\kappa\geq 2$,
$\lambda\geq 3$ and a minimal number of edges.}\label{5}
\end{figure}

{\bf Six vertices.} From the classification of graphs on six
vertices (see for instance~\cite{Speed}) we know that any such
2-vertex and 3-edge connected graph can be obtained by adding
edges to one of the four graphs shown in Figure~\ref{5}. By
Corollary~\ref{Corollary 1.8.} it suffices to check the formula
of the proposition for them.

\begin{figure}
$$\epsfbox{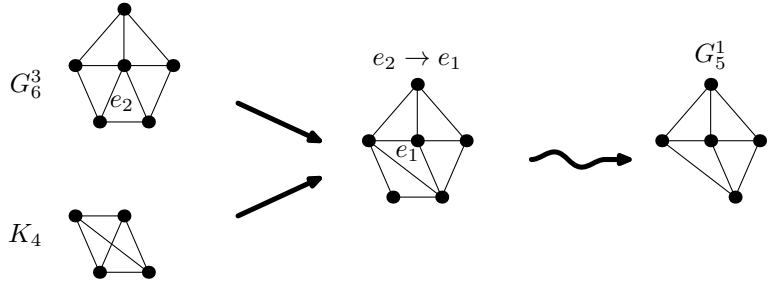}$$
\caption{Using Proposition~\ref{Proposition 1.9.} we get that
$\tau_2(G_6^3)=\tau_2(G_5^1)=1$.}\label{6}
\end{figure}

Note that $G_6^1$ and $G_6^2$ have 9 edges. They both have zero
$2$-TC ($9{-}12{+}3=0$). Indeed, in Subsection~\ref{A simple
examples} we mentioned that $B_2(G_6^1)$ has codimension 1 strata
with nontrivial linear fiber. As it is stated in~\cite{WW} the
graph $G_6^2$ underlies a tensegrity if and only if the six
points lie on a conic, which is also a codimension~1 condition.
Note that $G_6^2$ is the complete bipartite graph $K_{3,3}$. For
$G_6^3$ we proceed as follows. From Corollary~\ref{Corollary
1.7.} it follows that $\tau_2(G_6^3)\geq 1$. Then we use
Proposition~\ref{Proposition 1.9.} in the same way as in
Example~\ref{Example 1.10.} to show that
$$
\tau_2(G_6^3)=\tau_2(G_5^1)=1, \quad \hbox{and again} \quad
10-12+3=1,
$$
see Figure~\ref{6}.

It is easy to see that the same argument works to show that
$\tau_2(G_6^4)=1$.

\vspace{1mm}

\begin{figure}
$$\epsfbox{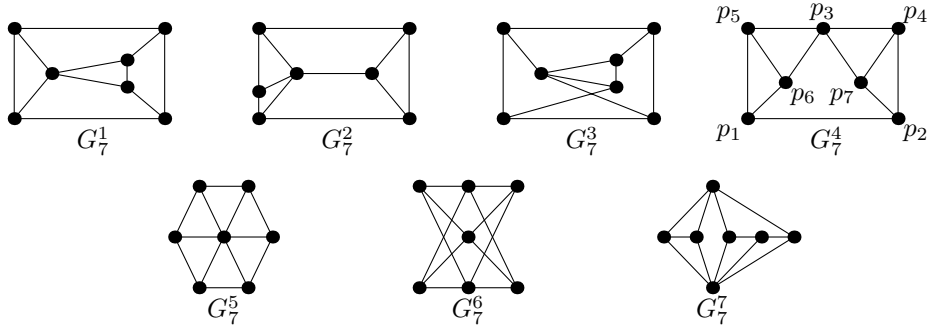}$$
\caption{The seven graphs with seven vertices, $\kappa\geq 2$,
$\lambda\geq 3$ and a minimal number of edges.}\label{7}
\end{figure}

{\bf Seven vertices.} From the classification of graphs with
seven vertices (see \cite{Speed}) we get that all 2-vertex and
3-edge connected graphs on seven vertices can be obtained by
adding edges to one of the seven graphs shown in Figure~\ref{7}.
By Corollary~\ref{Corollary 1.8.} it suffices again to check
these graphs.

To prove that $\tau_2(G_7^1)=\tau_2(G_7^2)=0$ we use
Corollary~\ref{Corollary 1.12.NEW} applied to $G_6^1$. Note that
the geometric conditions for $G_6^1$ to underlie a nonzero
tensegrity (see Subsection~5.1) allow to apply
Corollary~\ref{Corollary 1.12.NEW}. Similarly, we apply
Corollary~\ref{Corollary 1.12.NEW} to $G_6^2$ to conclude that
$\tau_2(G_7^3)=0$. By computations analogous to~\cite[Section
4]{Guz1} we find that $\tau_2(G_7^4)=0$. Indeed, one can show
that this graph underlies a nonzero tensegrity if and only if at
least one of the following codimension~1 conditions holds:
$$
p_1\triangledown p_2\triangledown p_3 =0, \quad
p_1\triangledown p_5\triangledown p_6 =0, \quad
p_2\triangledown p_4\triangledown p_7 =0, \quad
p_3\triangledown p_4\triangledown p_7 =0, \quad
p_3\triangledown p_5\triangledown p_6 =0.
$$

\begin{figure}
$$\epsfbox{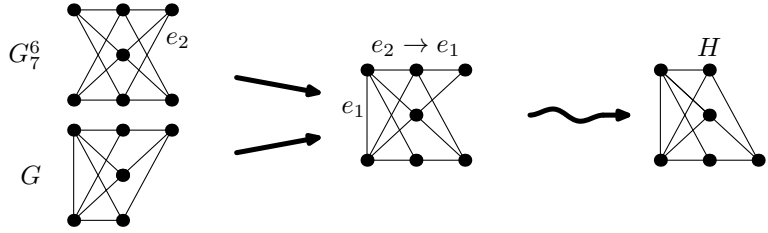}$$
\caption{Using Proposition~\ref{Proposition 1.9.} one sees that
$\tau_2(G_7^6)=1$.}\label{8}
\end{figure}

So the first four graphs with 11 edges have zero $2$-TC. The
other three have 12 edges. We apply Corollary~\ref{Corollary
1.12.NEW} to $G_6^3$ and $G_6^4$ to obtain that
$$
\tau_2(G_7^5)=1 \quad \hbox{and} \quad \tau_2(G_7^7)=1.
$$
To prove that the $2$-TC of $G_7^6=K_{3,4}$ is 1 we proceed as
follows. First, $\tau_2(G_7^6)\geq 1$ by Corollary~\ref{Corollary
1.7.}. Then we apply Proposition~\ref{Proposition 1.9.} as shown
in Figure~\ref{8}. The graph $G$ has 6 vertices and 10 edges and
thus we have $\tau_2(G)=1$. It is easy to check that for a
general position framework $G(P)$ with a nonzero self-stress, all
edges of $G(P)$ have nonzero stress. On the middle picture we get
a vertex of degree 2, so we reduce to the graph $H$ on the right.
Note that $H$ is isomorphic to $G$, so $\tau_2(H)=1$. Hence
$\tau_2(G_7^6)=1$ as well.

\qed

\vspace{2mm}




\subsection{Geometric conditions for realizability of plane
tensegrities for graphs with zero tensegrity
$2$-characteristic}\label{smallvert2}

Like in intersection theory of algebraic varieties, it often
happens that strata for a graph with negative $2$-TC are obtained
as intersections of closures of some strata of graphs with zero
$2$-TC. So the conditions for realizability of plane tensegrities
for graphs with zero $2$-TC are the most important. In this
subsection we give all the conditions for the zero $2$-TC graphs
with number of vertices not exceeding 8.

In practice one would like to construct a tensegrity without
struts or cables with zero tension. So it is natural to give the
following definition. We say that a graph $G$ is {\it visible} at
the framework $P$ if there exists a self-stress that is nonzero at
each edge of this framework.

\begin{remark}
Visibility restrictions remove many degenerate strata. For
instance if a zero $2$-TC graph $G$ has a complete subgraph on
vertices $v_1$, $v_2$, and $v_3$, then the codimension 1 stratum
defined by the condition: {\it the points $v_1$, $v_2$, and $v_3$
are on one line} does in many cases not contain visible
frameworks.
\end{remark}

Let us list the geometric conditions for the vertices of all
visible 2-vertex and 3-edge connected graphs with $n$ vertices
and zero $2$-TC for $n\le 8$. To find the geometric conditions we
essentially use the surgeries of Section~4, see
Propositions~\ref{Proposition 1.9.}, \ref{operation1}
and~\ref{operation2}.

\vspace{2mm}

In the next table we use besides the elementary also the following
two additional geometric conditions:

--- six points are on a conic;

--- for six points $p_1, \ldots, p_6$ the
lines $p_1p_2$, $p_3p_4$, and $p_5p_6$ have a common nonempty intersection.\\
As we have seen in Examples~\ref{6pconic} and~\ref{2_2_2} these
conditions are equivalent to geometric conditional systems.

\newpage

\begin{center}
\begin{tabular}{|p{35mm}|p{115mm}|}
\hline
{\bf Graph (6 vert.)}  &
{\bf Sufficient geometric conditions}\\
\hline \hline
$\begin{array}{l}\epsfbox{table.1}\end{array}$
& \vspace{-8mm} the lines $v_1v_2$, $v_3v_4$, and $v_5v_6$ have a
common nonempty intersection\\
\hline
$\begin{array}{l}\epsfbox{table.2}\end{array}$
& \vspace{-8mm} the six points $v_1$, $v_2$, $v_3$, $v_4$, $v_5$,
and
$v_6$ are on a conic\\
\hline
\end{tabular}
\end{center}

\begin{center}
\begin{tabular}{|p{35mm}|p{115mm}|}
\hline
{\bf Graph (7 vert.)}  &
{\bf Sufficient geometric conditions}\\
\hline \hline
$\begin{array}{l}\epsfbox{table.3}\end{array}$
& \vspace{-8mm} $v_1\triangledown v_2\triangledown v_3=0$\\
\hline
$\begin{array}{l}\epsfbox{table.4}\end{array}$
& \vspace{-8mm} the lines $v_1v_2$, $v_3v_4$, and $v_5v_6$
have a common nonempty intersection\\
\hline
$\begin{array}{l}\epsfbox{table.5}\end{array}$
& \vspace{-8mm} the lines $v_1v_2$, $v_3v_4$, and $v_5p$ where
$p=[v_2{,}v_6;v_3{,}v_7]$ have a common
nonempty intersection\\
\hline
$\begin{array}{l}\epsfbox{table.6}\end{array}$
& \vspace{-8mm} the six points $v_1$, $v_2$, $v_3$, $v_4$, $v_5$,
and
$p$, where $p=[v_1{,}v_6;v_3{,}v_7]$ are on a conic\\
\hline
\end{tabular}
\end{center}

\begin{center}
\begin{tabular}{|p{95mm}|p{55mm}|}
\hline
{\bf Graph (8 vert.)} &
{\bf Geometric conditions}\\
\hline \hline
$\begin{array}{l}\epsfbox{table.20}$
$\epsfbox{table.21}$
$\epsfbox{table.22}\\
\epsfbox{table.23}$ $\epsfbox{table.24}$
$\epsfbox{table.25}\\
\epsfbox{table.26}$ $\epsfbox{table.27}$
$\epsfbox{table.28}\end{array}$
& \vspace{-8mm} the lines $v_1v_2$, $v_3v_4$, and $v_5v_6$ have a
common nonempty intersection\\
\hline
$\begin{array}{l}\epsfbox{table.10}$ $\epsfbox{table.11}$
$\epsfbox{table.12}\end{array}$
& \vspace{-8mm} $v_1\triangledown v_2\triangledown v_3=0$\\
\hline
\end{tabular}
\end{center}

\begin{center}
\begin{tabular}{|p{95mm}|p{55mm}|}
\hline
{\bf Graph (8 vert.)} &
{\bf Geometric conditions}\\
\hline \hline
$\begin{array}{l}\epsfbox{table.30}\end{array}$
& \vspace{-8mm} the six points $v_1$, $v_2$, $v_3$, $v_4$, $v_5$,
and
$v_6$ are on a conic\\
\hline
$\begin{array}{l}\epsfbox{table.40}$
$\epsfbox{table.41}\end{array}$
& \vspace{-8mm} the lines $v_1v_2$, $v_3v_4$, and $v_5p$, where
$p=[v_2,v_6;v_3,v_7]$ have a
common nonempty intersection\\
\hline $\begin{array}{l}\epsfbox{table.50}\end{array}$
& \vspace{-8mm} the lines $v_1v_2$, $v_3v_4$, and $v_5p$, where
$p=[v_2,v_6;v_7,v_8]$ have a
common nonempty intersection\\
\hline $\begin{array}{l}\epsfbox{table.60}$
$\epsfbox{table.61}\end{array}$
& \vspace{-8mm} the six points $v_1$, $v_2$, $v_3$, $v_4$, $v_5$,
and
$p$, where $p=[v_1,v_6;v_3,v_7]$, are on a conic\\
\hline $\begin{array}{l} \epsfbox{table.70}\end{array}$
& \vspace{-8mm} the lines $v_1v_2$, $v_3p$, and $v_5q$, where
$p=[v_1,v_4;v_5,v_8]$ and $q=[v_2,v_6;v_3,v_7]$ have a
common nonempty intersection\\
\hline $\begin{array}{l} \epsfbox{table.80}\end{array}$
& \vspace{-8mm} the lines $v_5v_6$, $v_1p$, and $v_4q$, where
$p=[v_2,v_3;v_6,v_7]$ and $q=[v_2,v_3;v_6,v_8]$ have a
common nonempty intersection\\
\hline $\begin{array}{l} \epsfbox{table.90}\end{array}$
& \vspace{-8mm} the six points $v_1$, $v_2$, $v_4$, $v_6$, $p$,
and $q$, where $p=[v_2,v_3;v_6,v_7]$ and $q=[v_2,v_5;v_6,v_8]$,
are on
a conic\\
\hline $\begin{array}{l} \epsfbox{table.100}\end{array}$
& \vspace{-8mm} the six points $v_1$, $v_3$, $v_4$, $v_6$, $p$,
and $q$, where $p=[v_2,v_3;v_5,v_7]$ and $q=[v_5,v_7;v_6,v_8]$,
are on
a conic\\
\hline
 $\begin{array}{l} \epsfbox{table.110}\end{array}$
& \vspace{-8mm} the six points $v_1$, $v_2$, $v_3$, $v_5$, $p$,
and $q$, where $p=[v_1,v_6;v_3,v_7]$ and $q=[v_3,v_4;v_5,v_8]$,
are on
a conic\\
\hline
 $\begin{array}{l} \epsfbox{table.120}\end{array}$
& \vspace{-8mm} the six points $v_1$, $v_2$, $v_3$, $v_5$, $v_6$,
and $q$, where $p=[v_1,q;v_3,v_4]$ and $q=[v_5,v_7;v_4,v_8]$, are
on
a conic\\
\hline
\end{tabular}
\end{center}

\begin{center}
\begin{tabular}{|p{95mm}|p{55mm}|}
\hline
{\bf Graph (8 vert.)} &
{\bf Geometric conditions}\\
\hline
 $\begin{array}{l} \epsfbox{table.130}\end{array}$
& \vspace{-8mm} the six points $v_1$, $v_2$, $v_4$, $v_5$, $p$,
and $q$, where $p=[v_1,v_6;v_5,v_8]$ and $q=[p,v_7;v_2,v_3]$, are
on
a conic\\
\hline $\begin{array}{l} \epsfbox{table.140}$
$\epsfbox{table.141}\end{array}$
& \vspace{-8mm} the three points $[v_1,v_4;v_2,v_3]$,
$[v_1,v_5;v_2,v_6]$, and $[v_5,v_8;v_6,v_7]$ are on one
line\\
\hline
\end{tabular}
\end{center}

\begin{center}
\begin{tabular}{|p{35mm}|p{11.5cm}|}
\hline
{\bf Graph (8 vert.)} &
{\bf Sufficient geometric conditions}\\
\hline \hline
$\begin{array}{l}\epsfbox{table.150}\end{array}$
& \vspace{-8mm} the three points $[v_1,v_2;v_6,v_7]$,
$[v_1,p;v_6,v_8]$, and $[p,q;v_3,v_8]$,  where
$p=[v_2,v_4;v_5,v_8]$ and $q=[v_1,v_5;v_3,v_4]$, are on one line,
AND
the lines $p'v_2$, $q'v_3$, and $v_6v_7$ have a common nonempty
intersection, where $p'=[r',s';v_1,v_6]$, $q'=[r',s';v_6,v_8]$,
$r'=[v_1,v_4;v_2,v_5]$, and $s'=[v_3,v_4;v_5,v_8]$
\\
\hline
$\begin{array}{l}\epsfbox{table.160}\end{array}$
& \vspace{-8mm} the six points $v_1$, $v_4$, $v_7$, $v_8$, $p$,
and $q$, where $p=[r,s;v_3,v_4]$, $q=[r,s;v_5,v_8]$,
$r=[v_1,v_2;v_5,v_6]$, and $s=[v_2,v_3;v_6,v_7]$, are on a conic,
AND
the six points $v_1$, $v_2$, $v_6$, $v_7$, $p'$, and $q'$, where
$p'=[r',s';v_2,v_3]$, $q'=[r',s';v_5,v_6]$,
$r'=[v_1,v_4;v_5,v_8]$, and
$s'=[v_3,v_4;v_7,v_8]$, are on a conic\\
\hline $\begin{array}{l}\epsfbox{table.170}\end{array}$
& \vspace{-8mm} the six points $v_1$, $v_3$, $v_4$, $v_5$, $v_7$,
and $p$, where $p=[v_1,q;v_7,v_8]$, $q=[r,s;v_2,v_3]$,
$r=[v_3,v_6;v_7,v_8]$, and $s=[v_1,v_6;v_2,v_8]$, are on a conic,
AND
the six points $v_1$, $v_2$, $v_3$, $v_6$, $v_8$, and $p'$, where
$p'=[v_3,q';v_7,v_8]$, $q'=[r',s';v_1,v_4]$,
$r'=[v_1,v_5;v_7,v_8]$, and $s'=[v_3,v_5;v_4,v_7]$, are on a conic

\\
\hline
\end{tabular}
\end{center}

\begin{remark}
For the last three graphs in the table we have two distinct
equations. Nevertheless, the $2$-TC of the graphs are zero. This
is similar to the case of non-complete intersections in algebraic
geometry.
\end{remark}

\vspace{5mm}

\end{document}